\documentclass[11pt,a4paper]{article}
\usepackage{amsmath}
\usepackage{amssymb}
\usepackage{latexsym}

\newtheorem{theorem}{Theorem}
\newtheorem{corollary}{Corollary}
\newtheorem{proposition}{Proposition}

\newcommand{\prf}{{\em Proof}. }
\newcommand{\qed}{\hspace*{\fill}$\Box$}

\newcommand{\IR}{{\mathbb R}}
\newcommand{\IZ}{{\mathbb Z}}
\newcommand{\IC}{{\mathbb C}}

\newcommand{\kc}{{\cal C}}
\newcommand{\kk}{{\cal K}}
\newcommand{\kp}{{\cal P}}

\begin{document}

{\centerline{\Large{{\it Erratum:~}{\bf Compact hyperk\"ahler manifolds: basic
results}}}}

\bigskip

{\centerline{\large\bf Daniel Huybrechts}}

\bigskip

\bigskip

\bigskip

It was pointed out by D.\ Kaledin that the proof of Prop.\ 3.8 is
wrong. Actually, the proposition itself cannot be true as we
shall explain below. It was used to prove Cor.\ 3.10
and Thm.\ 3.11. The latter is the so called
projectivity criterion for hyperk\"ahler manifolds.
Here we will give a correct proof of these two
results. The main input is a recent theorem of J.-P.\ Demailly
and M.\ Paun \cite{DP}. All the other results of
the paper remain unaffected.

\bigskip

\bigskip

{\centerline{\bf The projectivity criterion for hyperk\"ahler manifolds }}
{\centerline{\bf as a consequence of the Demailly-Paun theorem}}

\bigskip

\noindent
{\bf 1.} Let $X$ be a compact hyperk\"ahler manifold of complex dimension
$2n$ and let ${\cal X}\to{\rm Def}(X)$ be the universal deformation of $X$.
For any cohomology class $\beta\in H^{4p}(X,\IR)$ let
$S_\beta\subset{\rm Def}(X)$ be the set of those $t\in{\rm Def}(X)$
for which $\beta$ is a cohomology class of type $(2p,2p)$ on
${\cal X}_t$. Then $S_\beta$ is a closed analytic subset of ${\rm Def}(X)$.
Let ${\rm A}\subset H^*(X,\IZ)$ be the set of all integral classes
$\beta\in H^{4p}(X,\IZ)$, $p=1,\ldots, n$ such that $S_\beta$
is a proper subset of ${\rm Def}(X)$. A point $t\in{\rm Def}(X)$ is called
{\it very general} if $t$ is in the complement of
$\bigcup_{\beta\in{\rm A}} S_\beta$ and if ${\cal X}_t$ does not admit any
analytic subsets of odd dimension.

\bigskip

{\it The set of very general points in ${\rm Def}(X)$ is dense.
If $t\in{\rm Def}(X)$ is a very general point and
$\beta\in H^{2p,2p}({\cal X}_t,\IZ)$ is an integral class of type
$(2p,2p)$ on ${\cal X}_t$, then $\beta$ is of type $(2p,2p)$ 
on any small deformation of ${\cal X}_t$. }

\bigskip

Here we use Fujiki's result \cite[Prop.\ 5.11]{Fujiki}
and the fact that $\bigcup_{\beta\in{\rm A}} S_\beta$
is a countable union of proper closed analytic subsets.
The latter says that
for any hyperk\"ahler metric on $X$ the general complex structure
compatible with it does not admit any odd-dimensional analytic subset. 

For classes $\beta$ which are of pure type $(2p,2p)$ on
any small deformation of $X$ the form of degree $2n-2p$
on $H^2(X,\IC)$ defined by $\alpha\mapsto \int_X\beta\alpha^{2(n-p)}$
can be expressed in terms of the Beauville-Bogomolov quadratic form
$q_X$ on $H^2(X,\IC)$. More precisely we have:

\bigskip

{\it If $\beta\in H^{4p}(X,\IC)$ is of type $(2p,2p)$ on all
small deformations of $X$, then there exists a constant
$c_\beta$ depending on $\beta$ such that for all $\alpha\in H^2(X,\IC)$
one has $\int_X\beta\alpha^{2(n-p)}=c_\beta q_X(\alpha)^{n-p}$.}

\bigskip

This is Thm.\ 5.12 in \cite{Huyhab}, which is a generalization of a
result of Fujiki and whose proof uses arguments of Bogomolov
(cf.\ \cite[1.11]{Huyinv}).

\bigskip

\bigskip

\noindent
{\bf 2.} We call the compact hyperk\"ahler manifold $X$ itself {\it very
general} if the point $0\in{\rm Def}(X)$ corresponding to it is a very
general point.
As a consequence of the above one obtains:

\bigskip

{\it If $X$ is a very general compact hyperk\"ahler manifold and $Y\subset
X$ is an irreducible analytic subset then its codimension is even, say
$2p$, and the cohomology class $[Y]\in H^{2p,2p}(X,\IZ)$ is of type
$(2p,2p)$ on all small deformations of $X$. In particular, there exists
a constant $c_{[Y]}$ such that
$\int_Y\alpha^{2(n-p)}=c_{[Y]}q_X(\alpha)^{n-p}$ for all
$\alpha\in H^2(X,\IC)$.}

\bigskip

We next quote the result of Demailly and Paun:

\begin{theorem}\cite{DP}---
Let $X$ be a compact K\"ahler manifold. Then the K\"ahler cone
$\kk_X$ of $X$ is a connected component of the set $\kp_X$ of all
classes $\alpha\in H^{1,1}(X,\IR)$ such that $\int_Y\alpha^d>0$
for any irreducible analytic subset $Y\subset X$ of dimension $d$.
\end{theorem}

Combining this with the above one obtains
a description of the K\"ahler cone of a very general
hyperk\"ahler manifold. Recall that the {\it positive
cone} $\kc_X\subset H^{1,1}(X,\IR)$ is 
the connected component of the open subset $\{\alpha\in
H^{1,1}(X,\IR)~|~q_X(\alpha)>0\}$ that contains the K\"ahler cone
$\kk_X$.

\begin{corollary}---
Let $X$ be a very general compact hyperk\"ahler manifold.
Then $\kk_X=\kc_X$. Moreover, the interior of the cone of
pseudo-effective classes coincides with $\kk_X=\kc_X$.
\end{corollary}

\prf One first shows that $\kc_X\subset\kp_X$. Since $\kk_X\subset\kp_X$
and $\kc_X$ is connected, it suffices to show that for any
$\alpha\in\kc_X$ and any irreducible analytic subset $Y\subset X$ of
codimension $2p$ one has $\int_Y\alpha^{2(n-p)}\ne0$. Since
$\int_Y\alpha^{2(n-p)}=c_{[Y]}q_X(\alpha)^{n-p}$ and $q_X(\alpha)>0$,
this follows from $c_{[Y]}\ne0$. The latter can be obtained
from the same equation applied to a K\"ahler class.
Thus,  $\kc_X\subset\kp_X$. Since $\kc_X$ is connected and contains
$\kk_X$, the Demailly-Paun theorem shows that $\kc_X=\kk_X$.

Clearly, every class in $\kk_X=\kc_X$ is in the interior of the
cone of pseudo-effective classes, which consists of the cone
of classes that can be represented by closed positive 
$(1,1)$-currents bounded from below by a K\"ahler form.
Conversely, if $\alpha$ can be represented
by a closed strictly positive $(1,1)$-current then $q_X(\alpha,\omega)=
c\int(\sigma\bar\sigma)^{n-1}\alpha\omega>0$ for any K\"ahler class
$\omega\in\kk_X=\kc_X$. Here, $c$ is a certain positive scalar and
$\sigma$ is a non-degenerate holomorphic two-form. Hence, $\alpha\in\kc_X$.
\qed

\bigskip

\noindent
{\bf 3.} Since the set of very general $t\in{\rm Def}(X)$ is dense,
any class $\alpha\in\kc_X\subset H^2(X,\IR)$ can be approximated by a sequence 
$\alpha_{t_i}\in H^2(X,\IR)$ such that $t_i$ is a sequence of very general
points converging to $0\in{\rm Def}(X)$ and $\alpha_{t_i}$ is a
K\"ahler class on ${\cal X}_{t_i}$. This is enough to conclude:

\begin{proposition}---
Let $X$ be a compact hyperk\"ahler manifold and let $\alpha\in\kc_X$
Then $\alpha$ is in the interior of the cone of pseudo-effective
classes, i.e.\ $\alpha$ can be represented by a closed
positive $(1,1)$-current which can be bounded from
below by a (small) K\"ahler form.
\end{proposition}

\prf Here one can copy the argument of Demailly \cite[Prop.\ 6.1]{Demailly}
that shows that the cone of pseudo-effective classes is closed.
Fix K\"ahler classes $\omega_t$ on ${\cal X}_t$ depending continously
on $t\in{\rm Def}(X)$. Then the mass of $\alpha_{t_i}$, which
is $\int_X\omega_{t_i}^{2n-1}\alpha_{t_i}$,
converges to $\int_X\omega_0^{2n-1}\alpha$. Hence the sequence
of forms $(\alpha_{t_i})$ is weakly bounded and thus weakly compact.
In particular, it contains a weakly convergent subsequence. As the
$\alpha_{t_i}$ are closed positive $(1,1)$-forms on ${\cal X}_{t_i}$,
the limit current is closed and positive of bidegree $(1,1)$ on $X$.
As $\kc_X$ is open, the class $\alpha$ must be in the interior
of the cone of all pseudo-effective classes.\qed

\bigskip

As a consequence we obtain the projectivity criterion for compact
hyperk\"ahler manifolds which was stated as Thm.\ 3.11 in \cite{Huyinv},
but the proof of which was seriously flawed as it used the wrong Prop.\
3.8.

\begin{theorem}---
Let $X$ be a compact hyperk\"ahler manifold. Then $X$ is projective if and
only if there exists a line bundle $L$ on $X$ with $q_X(c_1(L))>0$.
\end{theorem}

\prf If $X$ is projective, then there exists an ample line bundle $L$.
As $c_1(L)$ is then a K\"ahler class
and $\kk_X\subset\kc_X$, this yields $q_X(c_1(L))>0$.
Conversely, let us assume that there exists a line bundle $L$
with $q_X(c_1(L))>0$. This is equivalent to the existence of a line bundle
$L$ with  $c_1(L)\in\kc_X$. By the previous proposition $c_1(L)$
is in the interior of the cone of pseudo-effective classes.
Thus $c_1(L)$ can be represented by a closed positive $(1,1)$-current
which is bounded from below by a K\"ahler form.
Applying results of Bonavero \cite{Bonavero}
and Ji-Shiffman \cite{Shiffman} one obtains that $X$
is Moishezon and hence projective.\qed

\bigskip

\noindent
{\bf 4.} Let me indicate why Prop.\ 3.8 has to be false.
One way to see this is to use the fact that from dimension four on,
the birational K\"ahler cone can really be different
from the K\"ahler cone itself (something that does not happen for K3
surfaces). In this case one just picks a class that is positive
on the K\"ahler cone, but not on the entire birational K\"ahler cone.
If the original Prop.\ 3.8 were true, this class would be representable
by a closed positive current. But considered on any other birational compact
hyperk\"ahler manifold it would also be representable by a closed
positive current, which is in contradiction with the fact that there
is at least on birational compact hyperk\"ahler manifold where
the class is not positive on the K\"ahler cone. In this argument we
use that the quadratic form $q_X$ is compatible with
birational correspondences.

\bigskip

\noindent
{\bf 5.} Some of the results of the paper have been strengthened
in the meantime. In \cite{Huy2} we prove Thm.\ 4.6. without assuming 
the projectivity of the varieties. A better description of the
closure of the K\"ahler cone than the one given in Thm.\ 7.1
has also been found: A class is in the closure of the K\"ahler cone
if and only if it is non-negative on all rational curves (Prop.\ 3.1
in \cite{Huy2}).

\bigskip

\bigskip

\noindent
{\bf Acknowledgement:} First of all I wish to thank D.\ Kaledin for
pointing out the blunder in the proof of Prop.\ 3.8 and for believing
in the projectivity criterion and its consequences all the time. He was
always very patient to check various arguments over the last two years.
I am certainly most grateful to J.-P.\ Demailly and M.\ Paun for
their beautiful theorem. Thanks to J.-P.\ Demailly for comments
on the new proof.

\bigskip

\bigskip

\noindent
{\small Mathematisches Institut\\
Universit\"at zu K\"oln\\
Weyertal 86-90\\
50931 K\"oln, Germany\\
\texttt{huybrech@mi.uni-koeln.de}}


\bigskip

\bigskip

{\footnotesize
\begin{thebibliography}{mm}

\bibitem{Bonavero} Bonavero, L. \em In\'egalit\'es de Morse holomorphes
singuli\`eres. \em C.\ R.\ Acad.\ Sci.\ Paris 317 (1993),
1163-1166.

\bibitem{DP} Demailly, J.-P., Paun, M. \em Numerical
characterization of the K\"ahler cone of a compact K\"ahler manifold. \em
math.AG/0105176, (2001).

\bibitem{Demailly} Demailly, J.-P. \em Regularization of closed positive
currents and intersection theory. \em J.\ Alg.\ Geom.\ 1 (1992), 361-409

\bibitem{Fujiki} Fujiki, A. \em On the de Rham Cohomology Group of a Compact
K\"ahler Symplectic Manifold. \em Adv.\ Stud.\ Pure Math.\ 10 (1987),
105-165.

\bibitem{Huyinv} Huybrechts, D. \em Compact Hyperk\"ahler Manifolds:
Basic Results. \em Invent.\ math.\ 135 (1999), 63-113.

\bibitem{Huyhab} Huybrechts, D. \em Compact hyperk\"ahler manifolds. \em
Habilitation (1997).
http://www.mi.uni-koeln.de/\\
\~{}huybrech/artikel.htm/HKhabmod.ps

\bibitem{Huy2} Huybrechts, D. \em The K\"ahler cone of a compact
hyperk\"ahler manifold. \em math.AG/9909109

\bibitem{Shiffman} Ji, S., Shiffman, B. \em Properties
of compact complex manifolds carrying closed positive currents.
\em J.\ Geom.\ Anal.\ 3 (1993), 37-61.



\end{thebibliography}}
\end{document}